\documentclass[twoside,a4paper,12pt,centertags,reqno]{amsart} % 13pt displays smaller and better than 12pt
\usepackage{amsmath,amssymb,verbatim,vmargin}
\usepackage{color}
\usepackage{tikz}
\usepackage{tikz-cd}
\usetikzlibrary{matrix,calc}
\usepackage{hyperref}
\hypersetup{colorlinks,bookmarks=true,linktocpage=true,citecolor=blue,linkcolor=magenta}

\usepackage{eulervm}   % formulas round   

\allowdisplaybreaks % align break
\usepackage{color}
\usepackage[all]{xy} % diagrams

\usepackage{mathtools}
\usepackage{MnSymbol} %wedge righthalfcup

\theoremstyle{plain}
\newtheorem{thm}{Theorem}[section]

\theoremstyle{definition}

\newtheorem{rem}[thm]{Remark}

\usepackage{enumerate}
\usepackage{enumitem}

\newcommand{\bRn}{\mathbb{R}^n}

\newcommand{\loc}{\text{{\rm loc}}}

\newcommand{\pd}{\partial}

\newcommand{\bR}{{\mathbb R}}

\newcommand{\bN}{{\mathbb N}}

\newcommand{\cM}{{\mathcal M}}

\newcommand{\fH}{{\mathbf H}}

\def\barint_#1{\mathchoice
            {\mathop{\vrule width 6pt
height 3 pt depth -2.5pt
                    \kern -9.5pt
\intop \kern -4pt}\nolimits_{#1}}%
            {\mathop{\vrule width 5pt height
3 pt depth -2.6pt
                    \kern -6.5pt
\intop \kern -4pt}\nolimits_{#1}}%
            {\mathop{\vrule width 5pt height
3 pt depth -2.6pt
                    \kern -6pt
\intop \kern -4pt}\nolimits_{#1}}%
            {\mathop{\vrule width 5pt height
3 pt depth -2.6pt
          \kern -6pt \intop \kern -4pt}\nolimits_{#1}}}
          
           \def\bariint_#1{\mathchoice
            {\mathop{\vrule width 15pt
height 3 pt depth -2.5pt
                    \kern -15.8pt
\intop \kern -8pt\intop \kern -4pt}\nolimits_{#1}}%
            {\mathop{\vrule width 9pt height
3 pt depth -2.6pt
                    \kern -10.5pt
\intop \kern -8pt\intop \kern -4pt}\nolimits_{#1}}%
            {\mathop{\vrule width 9pt height
3 pt depth -2.6pt
                    \kern -10pt
\intop \kern -8pt\intop \kern -4pt}\nolimits_{#1}}%
            {\mathop{\vrule width 9pt height
3 pt depth -2.6pt
          \kern -8pt \intop \kern -10pt\intop \kern -4pt}
      \nolimits_{  #1}}}

\def\barintlim_#1{\mathchoice
            {\mathop{\vrule width 6pt
height 3 pt depth -2.5pt
                    \kern -8.8pt
\intop \kern -4pt}\limits_{#1}}%
            {\mathop{\vrule width 5pt height
3 pt depth -2.6pt
                    \kern -6.5pt
\intop \kern -4pt}\limits_{#1}}%
            {\mathop{\vrule width 5pt height
3 pt depth -2.6pt
                    \kern -6pt
\intop \kern -4pt}\limits_{#1}}%
            {\mathop{\vrule width 5pt height
3 pt depth -2.6pt
          \kern -6pt \intop \kern -4pt}\limits_{#1}}}
          
           \def\bariintlim_#1{\mathchoice
            {\mathop{\vrule width 15pt
height 3 pt depth -2.5pt
                    \kern -15.8pt
\intop \kern -8pt\intop \kern -4pt}\limits_{#1}}%
            {\mathop{\vrule width 9pt height
3 pt depth -2.6pt
                    \kern -10.5pt
\intop \kern -8pt\intop \kern -4pt}\limits_{#1}}%
            {\mathop{\vrule width 9pt height
3 pt depth -2.6pt
                    \kern -10pt
\intop \kern -8pt\intop \kern -4pt}\limits_{#1}}%
            {\mathop{\vrule width 9pt height
3 pt depth -2.6pt
          \kern -8pt \intop \kern -10pt\intop \kern -4pt}
      \limits_{  #1}}}
          
\renewcommand{\iint}{\int \kern -8pt\int}       

\numberwithin{equation}{section}
\makeatletter
\@namedef{subjclassname@2020}{\textup{2020} Mathematics Subject Classification}
\makeatother

\title{Maximal Regularity under Quadratic Estimates}

\author{Yi C. Huang} 
\address{School of Mathematical Sciences, Nanjing Normal University, Nanjing 210023, People's Republic of China}
\email{Yi.Huang.Analysis@gmail.com}
\urladdr{https://orcid.org/0000-0002-1297-7674}

\date{\today} 

\thanks{Research of the author is supported by the National NSF grant of China (no. 11801274).
The author would like to thank Li Liu (YZU) for calling his attention back to Weighted Maximal Regularity,
and Jian-Hua Chen (HNUST) for highlighting quadratic estimates in Control Theory.}

\begin{document}

\begin{abstract}
In this Short Note we complement the intriguing harmonic analytic perspective due to P. Auscher and A. Axelsson for the abstract evolution equations. 
This concerns a unified approach to temporally weighted estimates for the forward and backward maximal regularity operators in presence of quadratic estimates and functional calculi.
In particular we provide several invariance properties for the maximal regularity operators either in evolution form or in balayage form.
\end{abstract}

\maketitle

%\tableofcontents

\section{Introduction}

Assume that $-A$ is a densely defined, closed linear operator on a Hilbert space $\fH$, with domain $D(A)$ and generating a bounded analytic semigroup $\{e^{-zA}: |\arg z|<\delta\}$, $0<\delta<\pi/2$.
In this Note we study the autonomous evolutions on $\bR_+=(0,\infty)$:
\begin{equation} \label{e:Evol}
\begin{cases}
\pd_tu+Au=f,\\
 u(0)=0
\end{cases}
\end{equation}
and
\begin{equation} \label{e:Evol'}
\begin{cases}
\pd_tv-Av=f,\\
 v(\infty)=0.
\end{cases}
\end{equation}
The so-called maximal regularity for the evolution equations \eqref{e:Evol}-\eqref{e:Evol'} in the Hilbert space $\fH$
(also referred to as de Simon's theorem \cite{dSim64}) states as follows: 

\bigskip 

\begin{quote}
there is some constant $C=C(A,\fH)>0$ such that for all $f\in L^2(\bR_+;\fH)$,
\begin{equation} \label{e:QMR-C}
|||Au|||\leq C|||f|||\quad\text{ and }\quad |||Av|||\leq C|||f|||,
\end{equation}
where $|||\cdot|||=\|\cdot\|_{L^2(\bR_+;\fH)}$.
\end{quote}

\bigskip

\noindent Thus via \eqref{e:Evol} (respectively, \eqref{e:Evol'}) and \eqref{e:QMR-C}, the solution $u$ (respectively, $v$) lies in $L^2(\bR_+;D(A))\cap H^1(\bR_+;\fH)$.
The maximal regularity conveyed in \eqref{e:QMR-C} means that $\pd_tu$ and $Au$ (respectively, $\pd_tv$ and $Av$) enjoy the same regularity as the source $f$.

There are many machineries proving de Simon's maximal regularity theorem, typically via Fourier transform,
and various extensions to Banach spaces, typically in $L^q(\bR_+;L^p(\Omega))$, where $\Omega\subset\bRn$. 
For more comprehensive materials about maximal regularity estimates and applications to evolution equations,
see \cite{DenHiePru03, Are04, KunWei04, Mon09} by Denk-Hieber-Pr\"uss, Arendt, Kunstmann-Weis, and Monniaux. 

\section{Invariance properties for maximal regularity operators}

Our aim is to complement the harmonic analytic perspective due to Auscher and Axelsson \cite{AusAxe11p} 
in their unified approach to the $L^2(\bR_+,t^{\pm1}dt;\fH)$ estimates for the maximal regularity operators in presence of quadratic estimates. 
Weighted maximal regularity estimates are motivated by elliptic boundary value problems \cite{AusAxe11},
and are useful also in the study of abstract evolution equations, see \cite[Section 2]{AusAxe11p}.

Recall that the forward maximal regularity operator $\cM_+$ is defined by
\begin{equation} \label{e:M+}
\cM_+(f)(t)=\int_0^t Ae^{-(t-s)A}f(s)ds,
\end{equation}
while the backward maximal regularity operator $\cM_-$ is given by
\begin{equation} \label{e:M-}
\cM_-(f)(t)=\int_t^{\infty} Ae^{-(s-t)A}f(s)ds.
\end{equation}
In the notations $\cM_\pm$, we omit the dependence in $A$.  
The operators $\cM_\pm$ are associated to the evolution equations \eqref{e:Evol}-\eqref{e:Evol'} as for appropriate $f$, we have
$$Au=\cM_+(f)\quad\text{ and }\quad Av=-\cM_-(f).$$
Therefore, maximal regularity problems translate into boundedness of $\cM_\pm$,
which are typical examples of singular integral operators with operator-valued kernels.

Recall that $A$ is said to satisfy the quadratic estimate $(Q)_A$ if for all $h\in\fH$,
\begin{equation} \label{e:QE}
\big{|}\big{|}\big{|}sAe^{-sA}h\big{|}\big{|}\big{|}_{-1}\leq C\|h\|_{\fH},
\end{equation}
where $|||\cdot|||_{-1}=\|\cdot\|_{L^2(\bR_+,s^{-1}ds;\fH)}$.
One also refers to $\big{|}\big{|}\big{|}sAe^{-sA}h\big{|}\big{|}\big{|}_{-1}$ as square function norm of $h$,
a notion widely used for the abstract Hardy spaces in harmonic analysis and elliptic boundary value problems, 
see Auscher \cite{Aus07}, Auscher-McIntosh-Russ \cite{AusMcIRus08}, Hofmann-Mayboroda-McIntosh \cite{HofMayMcI11}, and Auscher-Stahlhut \cite{AusSta16}.

The main results of this Short Note can be summarized as below.

\begin{thm} \label{thm:main}
Let $N\in\bN_+=\{1,2,\cdots\}$ and $\cM_\pm$ be given as in \eqref{e:M+}-\eqref{e:M-}. We have

\bigskip

i) (Evolution Formulae) Suppose that $(Q)_A$ holds. For $f_0\in\fH$, then
\begin{equation} \label{e:MR+Evol}
\cM_+\left((sA)^Ne^{-sA}f_0\right)(t)=\frac{1}{N+1}(tA)^{N+1}e^{-tA}f_0,
\end{equation}
\begin{equation} \label{e:MR-Evol}
\cM_-\left(Ae^{-NsA}f_0\right)(t)=\frac{1}{N+1}Ae^{-NtA}f_0,
\end{equation}
and the two formulae hold respectively in $L^2(\bR_+,t^{-1}dt;\fH)$ and $L^2(\bR_+,tdt;\fH)$.

\bigskip

ii) (Balayage Formulae) Suppose that $(Q)_{A^*}$ holds. 
For $f\in L^2(\bR_+,t^{-1}dt;\fH)$, then
\begin{equation} \label{e:MR+Bala}
\int_0^\infty Ae^{-NtA}\cM_+(f)(t)dt=\frac{1}{N+1}\int_0^\infty Ae^{-NsA}f(s)ds;
\end{equation}
and for $f\in L^2(\bR_+,tdt;\fH)$, then
\begin{equation} \label{e:MR-Bala}
\int_0^\infty (tA)^{N}e^{-tA}\cM_-(f)(t)dt=\frac{1}{N+1}\int_0^\infty (sA)^{N+1}e^{-sA}f(s)ds;
\end{equation}
and both formulae hold weakly in $\fH$.

\bigskip

iii) (Endpoint Balayage Formulae) Suppose that $(Q)_{A^*}$ holds. 
For $f\in L^2(\bR_+,tdt;\fH)$, if in addition $\int_0^\infty e^{-sA}f(s)ds$ converges weakly in $\fH$, then
\begin{equation} \label{e:MR+BalaE}
\int_0^\infty e^{-NtA}\cM_+(f)(t)dt=\frac{1}{N+1}\int_0^\infty e^{-NsA}f(s)ds;
\end{equation}
and for $f\in L^2(\bR_+,t^{-1}dt;\fH)$, if in addition $\int_0^\infty Ae^{-sA}f(s)ds=0$, then
\begin{equation} \label{e:MR-BalaE}
\int_0^\infty t^{N-1}A^{N}e^{-tA}\cM_-(f)(t)dt=\frac1N\int_0^\infty s^{N}A^{N+1}e^{-sA}f(s)ds;
\end{equation}
and both formulae hold weakly in $\fH$.
\end{thm}

\begin{proof}
By formal computations plus Fubini theorem.
\end{proof}

The above formulae are best understood as invariant subspace properties of $\cM_\pm$ either in evolution form or in balayage form.
By \textit{evolution}, we mean an extension $E$ mapping boundary elements in $\fH$ to $L^2_{\loc}(\bR_+;\fH)$,
while by \textit{balayage} (or, \textit{sweeping}), we mean the dual mapping of $E$ that sends elements from $L^2_{\loc}(\bR_+;\fH)$ to boundary elements in $\fH$.
Both mappings are extremely useful in the abstract Hardy space theory, see \cite{AusMcIRus08, HofMayMcI11}.
For concrete balayage operators in connection with $\overline{\pd}$-equation and probability, see Amar-Bonami \cite{AmaBon79} and Labeye-Voisin \cite{LabVoi03}.

\section{Illustrations and further remarks}

Note that \eqref{e:QE} implies (see Albrecht-Duong-McIntosh \cite{AlbDuoMcI96}) that $\forall\,h\in\fH$,
\begin{equation} \label{e:QE'}
\big{|}\big{|}\big{|}\psi_1(sA)h\big{|}\big{|}\big{|}_{-1}:=\big{|}\big{|}\big{|}(sA)^Ne^{-sA}h\big{|}\big{|}\big{|}_{-1}\leq C\|h\|_{\fH}
\end{equation}
and 
\begin{equation} \label{e:QE''}
\big{|}\big{|}\big{|}\psi_2(sA)h\big{|}\big{|}\big{|}_{-1}:=\big{|}\big{|}\big{|}e^{-sA}\left(I-e^{-NsA}\right)h\big{|}\big{|}\big{|}_{-1}\leq C\|h\|_{\fH},
\end{equation}
where $N\in\bN_+=\{1,2,\cdots\}$.
The functions $\psi_1$ and $\psi_2$ in \eqref{e:QE'}-\eqref{e:QE''} decay at 0 and $\infty$. 

For convenience, we also introduce weighted spaces $\mathfrak{H}_{\alpha}=L^2(\bR_+,t^\alpha dt;\fH)$, $\alpha\in\bR$.
Recall the following weighted de Simon's theorem \cite{AusAxe11p} for \eqref{e:Evol}-\eqref{e:Evol'}: $\forall\,\alpha<1$,
$$\|Au\|_{\mathfrak{H}_{\alpha}}\leq C\|f\|_{\mathfrak{H}_{\alpha}}\quad\text{ and }\quad  \|Av\|_{\mathfrak{H}_{-\alpha}}\leq C\|f\|_{\mathfrak{H}_{-\alpha}},$$
which, as we mentioned before, are equivalent to
\begin{equation} \label{e:wdSim}
\|\cM_+(f)\|_{\mathfrak{H}_{\alpha}}\leq C\|f\|_{\mathfrak{H}_{\alpha}}\quad\text{ and }\quad  \|\cM_-(f)\|_{\mathfrak{H}_{-\alpha}}\leq C\|f\|_{\mathfrak{H}_{-\alpha}}.
\end{equation}
For the endpoint weights, we have the following beautiful reduction \cite{AusAxe11p}:
\begin{equation} \label{e:wdSim'}
\left\|\cM_+(f)-Ae^{-tA}\int_0^\infty e^{-sA}f(s)\right\|_{\mathfrak{H}_{1}}\leq C\|f\|_{\mathfrak{H}_{1}},
\end{equation}
\begin{equation} \label{e:wdSim''}
\left\|\cM_-(f)-e^{-tA}\int_0^\infty Ae^{-sA}f(s)\right\|_{\mathfrak{H}_{-1}}\leq C\|f\|_{\mathfrak{H}_{-1}}.
\end{equation}
Proofs of \eqref{e:wdSim}-\eqref{e:wdSim''} in \cite[Section 1]{AusAxe11p} are independent of quadratic estimates.

Now we illustrate the convergence issues in Theorem \ref{thm:main}:
applying \eqref{e:QE'} and \eqref{e:wdSim} we see that the formulae \eqref{e:MR+Evol}-\eqref{e:MR-Evol} hold respectively in 
$\mathfrak{H}_{-1}=L^2(\bR_+,t^{-1}dt;\fH)$ and $\mathfrak{H}_{1}=L^2(\bR_+,tdt;\fH)$;
using a duality argument, \eqref{e:QE'} for $A^*$, and \eqref{e:wdSim} we see that the formulae \eqref{e:MR+Bala}-\eqref{e:MR-Bala} hold weakly in $\fH$.
For \eqref{e:MR+BalaE}, note that \eqref{e:QE''} for $A^*$ implies
$$\begin{aligned}
&\int_0^\infty e^{-sA}f(s)ds-\int_0^\infty e^{-NsA}f(s)ds\\
&\qquad\qquad=\int_0^\infty e^{-sA}\left(I-e^{-(N-1)sA}\right)f(s)ds
\end{aligned}$$
converges weakly in $\fH$. 
Hence for the family of integrals $\int_0^\infty e^{-NsA}f(s)ds$, the weak convergence in $\fH$ for $N=1$ encompasses the convergence for other $N\geq2$.
Note that by \cite[Remark 1.7]{AusAxe11p} or \eqref{e:wdSim'}, the assumptions for \eqref{e:MR+BalaE} guarantee $\cM_+(f)\in L^2(\bR_+,tdt;\fH)$.
For \eqref{e:MR-BalaE}, by \cite[Proposition 1.6]{AusAxe11p} or \eqref{e:wdSim''}, the trace condition
$$\begin{aligned}
&\lim_{\tau\rightarrow0}\frac{1}{\tau}\int_\tau^{2\tau}\cM_-(f)(t)dt\\
&\qquad\qquad=\int_0^\infty Ae^{-sA}f(s)ds=0 \quad\text{in}\quad\fH
\end{aligned}$$
 plus $f\in L^2(\bR_+,t^{-1}dt;\fH)$ guarantee $\cM_-(f)\in L^2(\bR_+,t^{-1}dt;\fH)$. 
The weak convergence in $\fH$ for \eqref{e:MR-BalaE} follows upon using a duality argument together with \eqref{e:QE'} for $A^*$.

\begin{rem}
Independent of $(Q)_A$, \eqref{e:MR+Evol}-\eqref{e:MR-Evol} also hold strongly in $\fH$ for $\forall \,t>0$.
\end{rem}

\begin{rem}
We formulated the formulae \eqref{e:MR+Evol}-\eqref{e:MR-BalaE} for $\mathfrak{H}_{\pm1}$ under a parameter $N\in\bN_+$.
This uses \eqref{e:QE'}-\eqref{e:QE''}.
However, by McIntosh's bounded holomorphic functional calculus (see e.g. Haase \cite{Haas06}), 
more general holomorphic functions other than $\psi_1$ and $\psi_2$ can be considered.
We refrain ourselves from such generalizations,
and in particular from considering $(sA)^{|\alpha|} e^{-sA}$ and $\mathfrak{H}_{\alpha}$ for $0<|\alpha|<1$ (see Ros\'en\footnotemark
\footnotetext{Axelsson becomes Ros\'en in publications later than \cite{AusAxe11p}.} \cite{Ros12}).
\end{rem}

\begin{rem}
In view of the neat decomposition \eqref{e:wdSim'}, it was observed in \cite[Remark 1.7]{AusAxe11p} that under $(Q)_{A^*}$,
$\cM_+$ extends to a bounded operator on the subspace of those $f\in L^2(\bR_+,tdt;\fH)$ so that $\int_0^\infty e^{-sA}f(s)ds$ converges weakly in $\fH$.
On one hand, \textit{there is no simple description of this subspace} (see \cite[Remark 1.7]{AusAxe11p}),
if no further information on $\fH$ and $e^{-tA}$ is provided.
On the other hand, this puzzle of subspace description can be connected to the Carleson duality in harmonic analysis, 
see Hyt\"onen and Ros\'en \cite{HytRos12, HytRos18}.
Our formula \eqref{e:MR+BalaE} gives an indirect description of this subspace: it is an invariant subspace of $L^2(\bR_+,tdt;\fH)$ under $\cM_+$.
\end{rem}

\bigskip

\section*{\textbf{Compliance with ethical standards}}

\bigskip

\textbf{Conflict of interest} The author has no known competing financial interests or
personal relationships that could have appeared to influence this reported work.

\bigskip

\textbf{Availability of data and material} Not applicable.

\bigskip

\bibliographystyle{alpha}
 
\bibliography{Hua-MaxRegQuadratic}

\end{document}